\documentclass[journal]{IEEEtran}
\usepackage{amsmath,amssymb}
\usepackage{overpic}
\usepackage{rotating}

\usepackage{textcomp,xcolor}

\renewcommand{\IEEEproofindentspace}{0em}
\renewcommand{\IEEEQED}{\IEEEQEDoff}
\newcommand{\rev}[2]{\textcolor{black}{#2}} 
\newcommand{\SubC}{\mathbb{S}_C}
\newcommand{\SubB}{\mathbb{S}_B}

\DeclareMathOperator*{\argmax}{arg\rm{}max}

\DeclareMathOperator*{\tr}{tr}

\newcommand{\ind}{\gamma}
\newcommand{\Ind}{\boldsymbol{\gamma}}
\newcommand{\ba}{\mathbf{a}}
\newcommand{\bA}{\mathbf{A}}
\newcommand{\bB}{\mathbf{B}}
\newcommand{\bC}{\mathbf{C}}
\newcommand{\bI}{\mathbb{I}}
\newcommand{\bP}{\mathbb{P}}
\newcommand{\bQ}{\mathbf{Q}}
\newcommand{\bR}{\mathbf{R}}
\newcommand{\bU}{\mathbf{V}}
\newcommand{\bS}{\mathbf{\Sigma}}
\newcommand{\bTi}{\mathbf{\Psi}} 
\newcommand{\bT}{\mathbf{\Phi}} 

\newcommand{\bW}{\mathbf{W}}
\newcommand{\bx}{\mathbf{x}}
\newcommand{\by}{\mathbf{y}}
\newcommand{\bz}{\mathbf{z}}
\newcommand{\bu}{\mathbf{u}}
\newcommand{\be}{\mathbf{e}}
\newcommand{\reals}{\mathbb{R}}

\newtheorem{theorem}{Theorem}
\newtheorem{lemma}[theorem]{Lemma}
\usepackage{algorithm}
\usepackage{algorithmic}
\usepackage{tikz}

\usepackage[colorlinks=true, linktocpage=false]{hyperref}
\usepackage{graphicx}
\usepackage{subfig}
\graphicspath{{figures/}}

\title{\LARGE{Optimal Sensor and Actuator \rev{Placement }{Selection} using Balanced Model Reduction}}
\author{Krithika~Manohar,~
        J.~Nathan~Kutz,~\IEEEmembership{Member,~IEEE,}
        and~Steven~L.~Brunton,~\IEEEmembership{Member,~IEEE}\vspace{-.2in}
\thanks{The authors thank Bing Brunton, Eurika Kaiser, and Josh Proctor for valuable discussions.  SLB acknowledges support from AFOSR Grant FA9550-18-1-200. JNK acknowledges support from AFOSR Grant FA9550-19-1-0011.   
KM acknowledges support from NSF MSPRF Award No. 1803289.}%
\thanks{K. Manohar was with the Department
of Applied Mathematics, University of Washington, Seattle,
WA, 98195 USA e-mail: kmanohar@uw.edu.}
\thanks{J. N. Kutz and S. L. Brunton are with University of Washington.}
\vspace{-.1in}}
\date{August, 2018}
\begin{document}

\maketitle

\begin{abstract}
Optimal sensor and actuator selection is a central challenge in high-dimensional estimation and control. 
Nearly all subsequent control decisions are affected by these sensor/actuator locations, and optimal placement amounts to an intractable brute-force search among the combinatorial possibilities.  
In this work, we exploit balanced model reduction and greedy optimization to efficiently determine sensor and actuator selections that optimize observability and controllability. 
In particular, we determine locations that optimize scalar measures of observability and controllability via greedy matrix QR pivoting on the dominant modes of the direct and adjoint balancing transformations. 
Pivoting runtime scales linearly with the state dimension, making this method tractable for high-dimensional systems. 
The results are demonstrated on the linearized Ginzburg-Landau system, for which our algorithm approximates known optimal placements computed using costly gradient descent methods.
\end{abstract}

\begin{IEEEkeywords}
optimal control, balanced truncation, sensor selection, actuator selection, observability, controllability.\vspace{-.175in}
\end{IEEEkeywords}

\section{Introduction}
\begin{figure}
\begin{center}
\includegraphics[width=.45\textwidth]{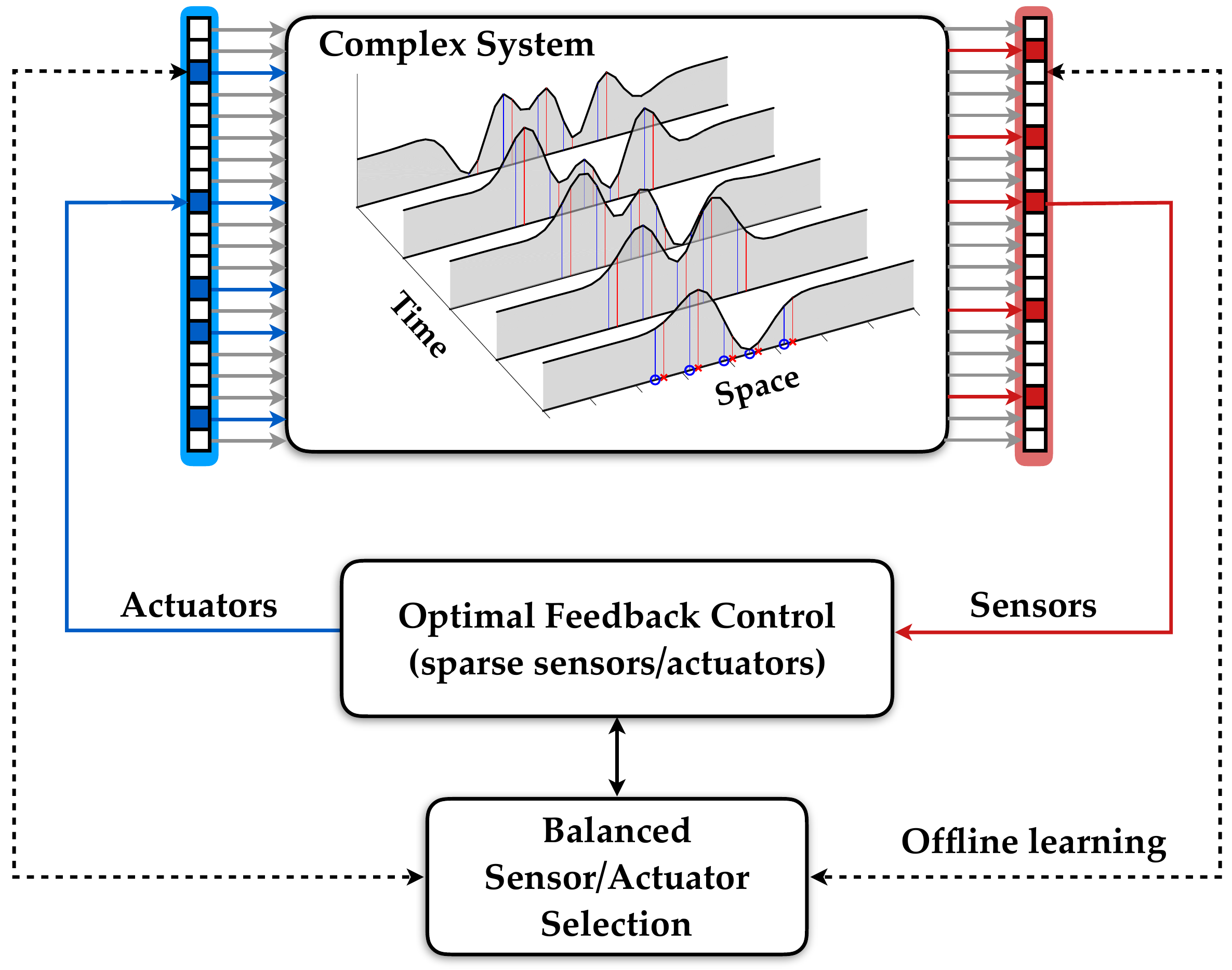}
\end{center}
\vspace{-.15in}
\caption{Schematic of balanced sensor and actuator selection for the optimal control of a high-dimensional system.}\label{Fig:Overview}
\vspace{-.15in}
\end{figure}

Optimizing the selection of sensors and actuators is one of the foremost challenges in feedback control~\cite{dp:book}. 
For high-dimensional systems it is impractical to monitor or actuate every state, hence a few sensors and actuators must be carefully positioned for effective estimation and control.
Determining optimal selections with respect to a desired objective is an NP-hard selection problem, and in general can only be solved by enumerating all possible configurations. 
This combinatorial growth in complexity is intractable; therefore, the placement of sensors and actuators are typically chosen according to heuristics and intuition.  
In this paper, we propose a greedy algorithm for sensor and actuator selection based on jointly maximizing observability and controllability in linear time-invariant systems. 
Our approach (see Fig.~\ref{Fig:Overview}) exploits low-rank transformations that balance the observability and controllability gramians to bypass the combinatorial search, enabling favorable scaling for high-dimensional systems.

To understand the challenges of sensor and actuator placement for estimation and control, we will first consider optimal sensor placement, which has mostly been used to reconstruct static signals.  
The primary challenge of sensor selection is that given $n$ possible locations and a budget of $r$ sensors, there are combinatorially many, $n\choose r$, configurations to evaluate in a brute-force search.  
Fortunately, there are heuristics that employ greedy selection of sensors based on maximizing mutual information~\cite{Krause2008jmlr} and information theoretic criteria~\cite{Paninski2005nc}. Another popular approach relaxes sensor selection to a weighted convex combination of possible sensors~\cite{Joshi2009ieee,chepuri2015continuous,Liu2016ieeetsp}, typically solved using semidefinite programming. Both heuristic approaches optimize submodular objective functions~\cite{summers2016submodularity}, which bound the distance between heuristic and optimal placement.  
Some objectives, such as those based on the quality of a Kalman filter, are not submodular~\cite{Zhang2017automatica}.  
Alternatively, sparsity-promoting optimization can be used to determine sensors and actuators~\cite{Lin2013ieeetac,Munz2014ieeetac,Zare2018arxiv}, although non-differentiability of sparsity promoting terms motivates other optimization techniques~\cite{Dhingra2014cdc}.

Even such heuristics cannot accommodate the large dimension of many physical models, such as in fluid dynamics. 
Fortunately, high-dimensional systems often evolve according to relatively few intrinsic degrees of freedom. 
Thus, it is possible to leverage dimensionality reduction to strategically select sensors. 
One  approach to place point sensors ~\cite{manohar2017data} computes the empirical interpolation points via EIM~\cite{Drmac2016siam} corresponding to the proper orthogonal decomposition (POD)~\cite{Berkooz1993pod} of data, to determine important locations in state space.

For systems with actuation, it is necessary to simultaneously consider the placement of sensors and actuators, since the most observable and most controllable subspaces are often different. 
Sensors and actuators for optimal feedback control are generally placed along the most observable and controllable directions, respectively~\cite{Chen:2011,Nestorovic2013mssp,hinson2014observability,Hemati2018actuator,summers2016submodularity}, using objective functions based on the associated observability or controllability gramians. 
Standard metrics for evaluating a certain sensor/actuator configuration include the $H_2$ norm~\cite{morris2015h2,Chen:2011}, a measure of the average impulse response, and the $H_\infty$ norm to measure the worst case performance. 
A chief drawback is the need to recompute the controller with each new configuration of sensors and actuators given by either the gradient minimization computation or brute-force searches. 
Moreover, these methods do not exploit the state-of-the-art in model reduction to optimize sensor and actuator placement.

\noindent\textbf{Contribution.} This work develops a scalable sensor and actuator selection algorithm based on balanced truncation~\cite{Moore1981ieeetac}, in which modes are hierarchically ordered by their observability and controllability. 
We use empirical interpolation of the low-rank balanced representation to find maximally observable and controllable states.
The resulting locations correspond to near-optimal point sensor and actuator configurations. 
The quality of our optimized configurations are evaluated using the $H_2$ norm of the resulting system, which is an average measure of its output energy. \rev{}{The closed loop $H_2$ norm is more relevant than open loop metrics for control performance, given a specific $H_2$ cost function.  Our approach, when used to optimize the open loop $H_2$ norm, is agnostic to the specific choice of controller weight matrices, and instead maximizes the input--output energy of the reduced order model. We also show that it is possible to apply our framework to closed loop systems, demonstrating near optimal sensor and actuator selection in comparison with more expensive iterative closed loop $H_2$ optimization. }
The runtime scales linearly with the number of state variables, after a one-time offline computation of the balancing transformation, which is less expensive than iterative alternatives. 
The resulting sensor and actuator configurations reproduce known optimal locations at a fraction of the cost associated with competing gradient descent methods.



\section{Problem Setup}

Consider the following \rev{stable}{} linear time-invariant system with a given state-space realization
\begin{subequations}\label{eqn:ss}
\begin{align}
	&\dot\bx = \bA\bx + \bB\bu  & \bx\in\reals^n,\bu\in\reals^q~ \\ 
	&\by = \bC\bx, &\by \in\reals^p,
\end{align}
\end{subequations}
\rev{with $q$ inputs and $p$ outputs. Sensors and actuators will be constrained to be sparse, localized points in state-space.  
Thus, the sensing and actuation matrices $\bC\in\reals^{p\times n}$ and $\bB\in\reals^{n\times q}$ must be structured in the following way}
{with large state dimension, i.e., $n\gg 1$. It is assumed that the system is stable, and $\bB$ and $\bC$ are linear {\em actuation} and {\em measurement} operators that make the system observable and controllable. Our objective is to choose a minimal subset of these sensors and actuators to obtain a system that is most jointly controllable and observable. For illustration we begin with $\bB=\bC=\mathbb{I}$ which correspond to pointwise sensing and actuation, but in general the subset selection can be adapted for arbitrary $\bB$ and $\bC$.
This subset selection corresponds to multiplying inputs and outputs by the selection matrices 
\begin{subequations}
\begin{align}
\SubC &= \begin{bmatrix} 
\be_{\ind_1} &
\be_{\ind_2} &
\dots &
\be_{\ind_r}
\end{bmatrix}^T \\
\SubB &= \begin{bmatrix} 
\be_{\beta_1} &
\be_{\beta_2} &
\dots &
\be_{\beta_r}
\end{bmatrix}.
\end{align}
\end{subequations}}%
Here $\be_j$ are the canonical basis vectors for $\reals^n$ with a unit entry at the \rev{}{selected} index $j$ and zeros elsewhere, where \rev{}{$\Ind = \{\ind_1,\dots,\ind_r\} \subset \{1,\dots,p\}$ denotes the index set of sensor locations with $\mbox{card}(\Ind)=r$.} Similarly, actuator selection indices are given by $\boldsymbol{\beta} = \{\beta_1,\dots,\beta_r\}$. 
\rev{}{The new measurement and actuation operators are $\hat\bC=\SubC\bC$ and $\hat\bB=\bB\SubB$ respectively. 
In the special case $\bB=\bC=\mathbb{I}$, the new operators $\hat\bC = \SubC\mathbb{I}$ and $\bB_\star = \mathbb{I}\SubB$ select subsets of state inputs and outputs,}
and the output would consist of $r$ components of $\bx$
\begin{equation}
\by = \hat\bC\bx = [x_{\ind_1} ~x_{\ind_2}~ \dots~ x_{\ind_r}]^T.
\end{equation}

\rev{A system with full actuation and sensing corresponds to $\bB=\bC=\mathbb{I}$.}{} Problem statement: {\em What are the best $r$-subsets of a given set of $p$ sensors and $q$ actuators, where $r\ll n$? }

To answer this question, we first quantify the degree of observability and controllability for a given set of sensors and actuators, i.e. for a given choice of $\bC$ and $\bB$. 
\rev{The matrices $\bC$ and $\bB$ consist of selected rows and columns of the $n\times n$ identity matrix.}{}  Optimizing over these directly involves a combinatorial search, and thus a heuristic approach is necessary for high-dimensional systems.

\subsection{Observability and controllability}

\begin{figure}
\begin{center}
\begin{overpic}[width=.475\textwidth]{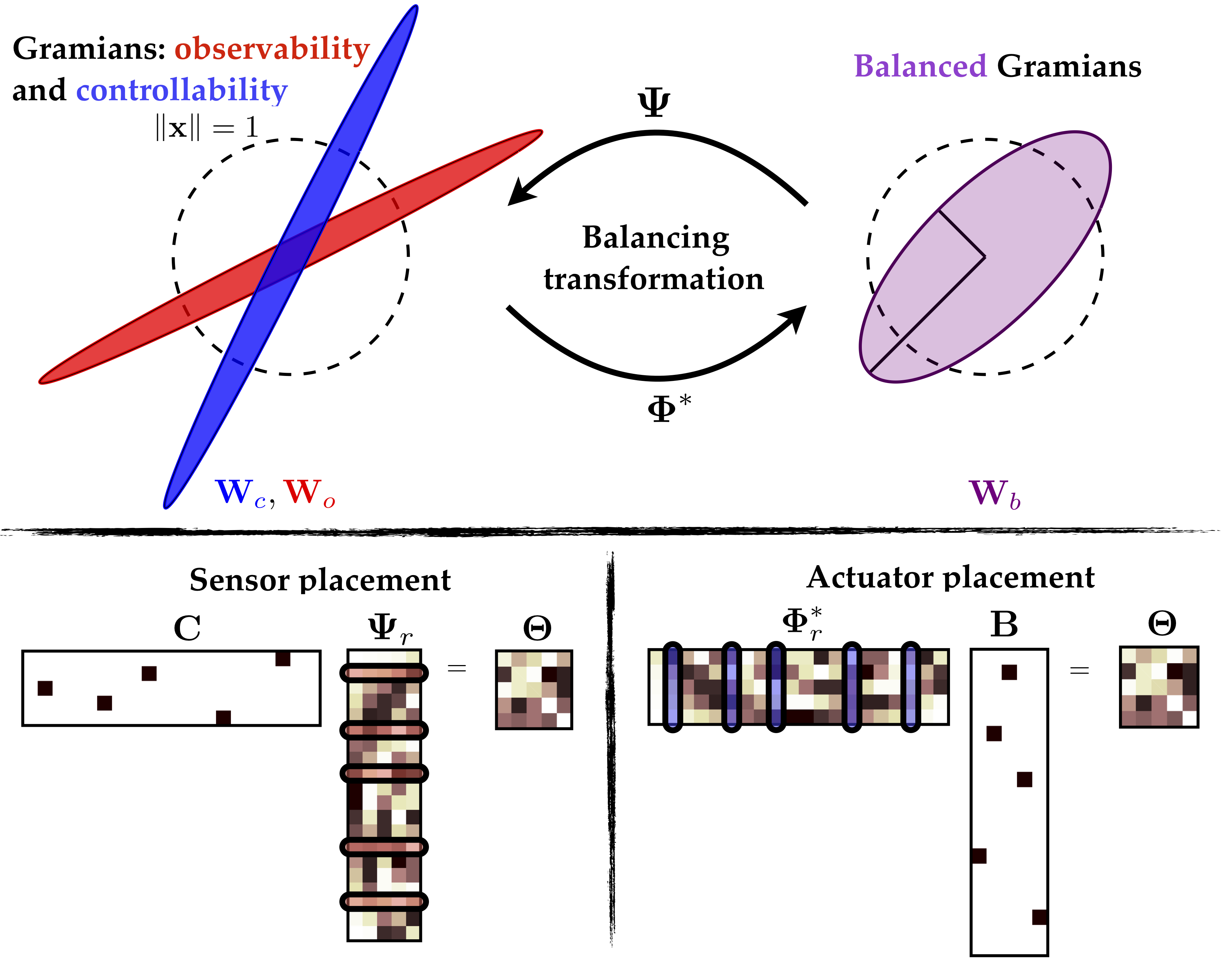}
\put(15.5,27.1) {$\hat{}$}
\put(81.5,27.5) {$\hat{}$}
\end{overpic}
\end{center}
\vspace{-.15in}
\caption{(top) Illustration of the balancing transformation for gramians.  The reachable set $\mathcal{E}_c$ with unit control input is shown in blue\rev{, given by $\bW_c^{1/2}\bx$ for $\|\bx\|=1$}{}.  The corresponding observable set is shown in red.  Under the balancing transformation $\bTi$, the gramians are equal, shown in purple.  (bottom) Sensor and actuator selection  based on balancing transformation.}\label{Fig:BTSenseAct}
\vspace{-.1in}
\end{figure}

The degrees of observability and controllability for the state-space system~\eqref{eqn:ss} are quantified by the observability gramian $\bW_o$ and controllability gramian $\bW_c$
\begin{equation}
\bW_o = \int_0^\infty e^{\bA^* t} \bC^* \bC e^{\bA t} dt, ~ 
\bW_c = \int_0^\infty e^{\bA t} \bB \bB^* e^{\bA^* t} dt,
\end{equation}
\rev{
Explicitly, the maximal energy output from a given initial condition $\bx_0$ is quantified by the observability gramian 
\begin{align}\label{eqn:max_out}
\max \|\by\|^2 &= \int_0^\infty \by(t)^*\by(t)dt \nonumber\\
&= \int_0 (\bC e^{\bA t} \bx_0)^*\bC e^{\bA t} \bx_0dt \nonumber \\
&= \bx_0^*\bW_o\bx_0.
\end{align}
Likewise, the minimal energy required to steer a system to a given state is defined by the inverse controllability gramian
\begin{equation}\label{eqn:min_in}
\min \|\bu\|_2 = \bx_0^*\bW_c^{-1}\bx_0.
\end{equation}
The gramians define ellipsoids that vary with the sensing matrix $\bC$ and actuation matrix $\bB$ (see Fig.~\ref{Fig:BTSenseAct}, top left).}{which may be visualized as controllable and observable \textit{ellipsoids} (Fig.~\ref{Fig:BTSenseAct}). These depend on the actuation and measurement operators, which consist of all states reachable from a bounded initial state
\begin{equation}
\mathcal{E}_c = \{\bW_c^{1/2}\bx \mid \|\bx\|_2 \le 1\},
\end{equation}
and all states that may be observed 
\begin{equation}
\mathcal{E}_o = \{\bW_o^{1/2}\bx \mid \|\bx\|_2 \le 1\}.
\end{equation}
}
Because the gramians depend on $\bB$ and $\bC$, they are often used to evaluate the observability/controllability of a given sensor and actuator placement. 
One important evaluation metric is the $H_2$ norm of a system. \rev{It measures the $L_2$ norm  or root mean square of the impulse response of a system, sometimes called its {\em output energy}. Given an impulsive input $u_j = \delta(0)$, the output response in component $i$ is ${y_{ij}(t)=\bC_i e^{\bA t}\bB_j}$. Therefore, the $H_2$ norm of \eqref{eqn:ss} is }{It measures the average output gain over all frequencies of the input, or the {\em output energy}. For the state-space system~\eqref{eqn:ss} with transfer function $G(s) = \bC(s\mathbb{I}-\bA)^{-1}\bB$, it is given by
\begin{equation}
\|G\|_2^2 = \frac{1}{4\pi^2}\int_0^\infty \tr(G(j\omega)^*G(j\omega)) d\omega.
\end{equation}
By the Plancherel theorem, it is also defined in the time domain by the impulse response ${y_{ij}(t)=\bC_i e^{\bA t}\bB_j}$ - the output in component $i$ given an impulse in input $j$, }
\begin{subequations}
\begin{align}\label{eqn:H2normB}
\|G\|_2^2 &=\int_0^\infty \tr(\bC e^{\bA t} \bB \bB^*e^{\bA^* t} \bC^*) dt = \tr(\bC\bW_c\bC^*)\\
 &= \int_0^\infty\tr(\bB^*  e^{\bA^* t} \bC^* \bC e^{\bA t}\bB) dt = \tr(\bB^*\bW_o\bB)
\end{align}
\end{subequations}
\rev{
In the content of sensor and actuator placement, the $H_2$ norm provides an objective function to be minimized over different choices of sensors and actuators, $\bC$ and $\bB$. 
In optimal control, it is desirable that the output energy for given impulse is made as small as possible, since this determines the ability of the controller to quickly stabilize a system upon excitement.
\subsection{A related metric}
A closely related metric to the $H_2$ norm is given by}{ which explicitly relate each gramian to {\em both} $\bB$ and $\bC$. A related alternative to the average output energy metric is given by the volumetric measure, the log determinant, denoted
\begin{equation}
\log|\bC\bW_c\bC^*|,\quad \log|\bB^*\bW_o\bB|,
\end{equation}
which are the logarithms of the geometric mean of the axes of the ellipsoid skewed by $\bB$ or $\bC$, by comparison the trace is the arithmetic mean. This metric is introduced by Summers et al~\cite{summers2016submodularity} to place actuators using a greedy optimization scheme for the submodular objective function
\begin{equation}
\bB_\star = \argmax_{\bB} \log|\bC\bW_c\bC^*|.
\end{equation}
For $H_2$ optimal control it is desirable to minimize the average gain from stochastic disturbance $\hat{w}$ to control output $\hat{z}(s) = \hat{G}(s)\hat{w}(s)$, namely, minimizing $\|\hat{G}\|_2$. Several strategies seek to build the controller and choose actuators simultaneously, using expensive gradient optimization schemes. The drawback of such closed loop metrics is having to recompute the gramians - an $O(n^3)$ operation - for every iteration that selects the next best actuator. This cubic scaling may be intractable for high-dimensional systems with large $n$. 
}

\rev{
Summers et al~\cite{summers2016submodularity} show that maximizing this log determinant is an effective proxy for $H_2$ norm minimization for choosing optimal actuators. We propose using this as an evaluation metric because it facilitates the efficient, greedy determinant maximization scheme detailed in Section~\ref{Sec:GreedyQR} which is solved using a greedy optimization scheme. 
\subsection{Our approach}
Our solution is to decouple sensor/actuator optimization from the gramian calculation, which is performed a single time. 
We do this by optimizing {\em sensors} using the full controllability gramian for $\bB=\mathbb{I}$, holding the actuators fixed
\begin{equation}
\bC_\star = \argmax_{\bC} \log|\bC\bW_c\bC^*|,
\end{equation}
omitting the dependence of $\bW_c$ on $\bB$ so that $\bW_c$ need not be recomputed at each step of the optimization. 
Likewise, we optimize {\em actuators} using the observability gramian for $\bC=\mathbb{I}$, holding the sensors fixed 
\begin{equation}
\bB_\star = \argmax_{\bB} \log|\bB^*\bW_o\bB|,
\end{equation}
in lieu of the trace metric~\eqref{eqn:H2normB}.
}{There are cases where optimizing sensors and actuators using the closed loop $H_2$ norm is more relevant for control~\cite{morris2015h2,Chen:2011}. 
By contrast, our approach reverses the strategy by instead starting from a maximally actuated and sensed optimal controller, then seeks a subset of these sensors/actuators to preserve (maximize) the geometric control measure, namely
\begin{subequations}\label{eqn:logdet_full}
\begin{align}
{\SubC}_\star &= \argmax_{\SubC} \log|\SubC\bC\bW_c\bC^*\SubC^*|, \\
{\SubB}_\star &= \argmax_{\SubB} \log|\SubB^T\bB^*\bW_o\bB\SubB|.
\end{align}
\end{subequations}%
} 
 Now, the gramians no longer depend on the optimization variable and need only be computed once, and both objectives are still fundamentally linked to the $H_2$ norm of the system. 
Critically, we will extract the dominant controllable and observable subspaces from a {\em balanced} coordinate transformation of the gramians.

\section{Balanced Model Reduction}
Many systems of interest are exceedingly high dimensional, making them difficult to characterize and limiting controller robustness due to significant computational time-delays.  
However, even if the ambient dimension is large, there may still be a few dominant coherent structures that characterize the system.   
Thus, significant effort has gone into obtaining efficient reduced-order models that capture the most relevant mechanisms for use in real-time feedback control~\cite{dp:book}.  

\rev{
The primary goal of balanced model reduction is to find a coordinate transformation $\bx\approx\bTi_r\ba_r$ giving rise to a related system $(\bA_r,\bB_r,\bC_r)$ with similar input--output characteristics, in terms of a state $\ba_r\in\mathbb{R}^r$ in a rank$-r$ basis $\bTi_r\in\mathbb{R}^{n\times r}$ with many fewer degrees of freedom, $r\ll n$. Instead of ordering modes based on energy, as in POD~\cite{Berkooz1993pod}, balanced model reduction identifies a coordinate transformation $\bx = \bTi \ba$ that hierarchically orders modes to capture the input--output characteristics of the system, as quantified by the controllability and observability gramians.}{ The goal of balanced model reduction is to find a transformation $\mathbf{T}$ from state-space (leaving inputs and outputs unchanged), $\begin{bmatrix}\begin{array}{c|c} \bA & \bB \\ \hline \bC & 0\end{array}\end{bmatrix}$ to $\begin{bmatrix}\begin{array}{c|c} \mathbf{T}\bA\mathbf{T}^{-1} & \mathbf{T}\bB \\ \hline \bC\mathbf{T}^{-1} & 0\end{array}\end{bmatrix}$, such that the transformed coordinates $\ba=\mathbf{T}^{-1}\bx$ are hierarchically ordered by their joint observability and controllability. This permits an $r$-dimensional representation made possible by truncating the $n-r$ least observable and controllable states.}
 
The seminal work of Moore in 1981~\cite{Moore1981ieeetac} showed it is possible to compute this coordinate system $\bTi$ where the controllability and observability gramians are equal and diagonal, denoted by the balanced model
\begin{align}
&\dot{\ba} = \bT^*\bA\bTi\ba + \bT^* \bB\bu  & {\ba}\in\reals^n,\bu\in\reals^q\nonumber \\ 
&\by = \bC\bTi{\ba}. &\by \in\reals^p
\end{align}
Here $\mathbf{T}^{-1}\triangleq\bTi$ are \emph{direct} modes and $\mathbf{T}\triangleq\bT^*$, the \emph{adjoint} modes.  
\rev{Note that $\bu$ and $\by$ are the same as in \eqref{eqn:ss}, even though the system state has been  reduced.}{}
The balanced state $\ba$ is then \emph{truncated}, keeping only the first $r\ll n$ most jointly controllable and observable states in $\ba_r$, so that $\bx\approx \bTi_r\ba_r$. This results in the \emph{balanced truncation} model~\cite{Moore1981ieeetac} $G_r=\begin{bmatrix}\begin{array}{c|c} \bT_r^* \bA\bTi_r & \bT_r^* \bB \\ \hline \bC\bTi_r & 0\end{array}\end{bmatrix}$. 
Since gramians depend on the particular choice of coordinate system, they will transform under a change of coordinates. 
The controllability and observability gramians for the balanced truncated system are
%
\begin{align}
\tilde\bW_c = \bT^*\bW_c\bT, \quad
\tilde\bW_o =  \bTi^*\bW_o\bTi.\label{Eq:TransWO}
\end{align}
%
The coordinate transformation $\bTi$ that makes the controllability and observability gramians equal and diagonal,
\begin{align}\label{Eq:BalancedgramiansHSV}
\tilde\bW_c = \tilde\bW_o = \bS,
\end{align}
is given by the matrix of eigenvectors of the product of the gramians $\bW_c\bW_o$ in the original coordinates:
\begin{align}
\hspace{-.1in}\tilde\bW_c\tilde\bW_o =\bT^*\bW_c\bW_o\bTi = \bS^2
~\Longrightarrow~  \bW_c\bW_o\bTi = \bTi\bS^2.\hspace{-.05in}\label{Eq:GRAMEIG}
\end{align}
\rev{}{The resulting balanced system is quantifiably close to the original system in the $H_\infty$ norm in terms of the {\em Hankel singular values} or diagonal entries of $\mathbf{\Sigma}$
\begin{equation}
\|G- G_r\|_\infty \le 2\sum_{i=r+1}^n \sigma_i.
\end{equation}
}%
In practice, computing the gramians $\bW_c$ and $\bW_o$ and the eigendecomposition of the product $\bW_c\bW_o$ in \eqref{Eq:GRAMEIG} may be prohibitively expensive for high-dimensional systems.  
Instead, the balancing transformation may be approximated with data from impulse responses of the direct and adjoint  systems, utilizing the singular value decomposition for efficient extraction of the relevant subspaces.  
The method of empirical gramians is quite efficient and is widely used~\cite{Moore1981ieeetac,Lall2002ijrnc,Willcox2002aiaaj,Rowley2005ijbc}.  
Moore's approach computes the entire $n\times n$ balancing transformation, which is not suitable for exceedingly high-dimensional systems.  
In 2002, Willcox and Peraire~\cite{Willcox2002aiaaj} generalized the method to high-dimensional systems, introducing a variant based on the rank-$r$ decompositions of $\bW_c$ and $\bW_o$ obtained from snapshots of direct and adjoint simulations.  
It is then possible to compute the eigendecomposition of $\bW_c\bW_o$ using efficient eigenvalue solvers.
This approach requires as many adjoint impulse-response simulations as the number of output equations, which may be prohibitively large for full-state measurements.  
In 2005, Rowley~\cite{Rowley2005ijbc} addressed this issue by introducing output projection, which limits the number of adjoint simulations to the number of relevant POD modes in the data. It is particularly advantageous to use these data-driven methods or low-rank alternating direction methods~\cite{li2002low} to approximate the gramians when there are fewer than full measurements and actuation of the state.

\section{\hspace{-.04in}Sensor \& actuator optimization via QR pivoting}\label{Sec:GreedyQR}

We now describe an efficient matrix pivoting algorithm to optimize the log determinant over the  choices of sensors and actuators.\rev{, described below and in algorithm~\ref{alg:qrpivot}.}{} 
The representation of the gramians in balanced truncation coordinates plays a crucial role. 
\subsection{Matrix volume objective}
Recall the goal of optimizing a set of \rev{$p$ point sensors and $q$ actuators out of $n$ possible choices.
Here we make the additional assumption that $p\ge r$ and $q\ge r$, where $r$ is the number of balanced modes required to faithfully approximate the full-order model. 
We begin with the balanced truncation ~\eqref{eqn:balred} of a state-space system~\eqref{eqn:ss}, and assume that the observable and controllable subspaces are well characterized by $r$ direct and adjoint modes $\bTi_r$ and $\bT_r$. }{$r$ sensors and actuators out of a fixed set $p$ and $q$ possible choices. The budget $r$ determines the balancing rank truncation, which necessarily must be less than both $p$ and $q$. Our sensor-actuator selection can be regarding as interpolating this rank-$r$ representation, that is, choosing locations or interpolation points that are heavily weighted in the dominant $r$ balanced modes. }

Summers et al~\cite{summers2016submodularity} show that it suffices to only consider controllable or observable subspaces for selecting sensors and actuators using the log determinant objective. Thus, 
we can substitute rank-$r$ balanced approximation of the gramians, $\hat\bW_c$ and $\hat\bW_o$, into the log determinant objective
\rev{}{
\begin{align}
\bC_\star &\approx \argmax_{\SubC} \log |\SubC\bC\bTi_r\bS_r\bTi_r^*\bC^T\SubC^T|  \nonumber \\
&= \argmax_{\SubC} |\SubC\bC\bTi_r|^2 \cdot |\bS_r| \nonumber \\
&= \argmax_{\SubC} | \SubC\bC\bTi_r|.
\end{align}}
This result follows from the monotonicity of logarithms and the product property of determinants, then omitting the term that is independent of the sensors, $\det\bS_r$. Likewise, in the actuator case, the objective $\argmax_{\SubB} \log|\hat\bB^T\hat{\bW}_o\hat\bB|$ simplifies 
\begin{equation}
\bB_\star = \argmax_{\SubB} |\bT^*\bB\SubB|.
\end{equation}
Consider for now the case of sensor placement.
The absolute determinant is a measure of matrix {\em volume}, and \rev{}{$\SubC$} is a row selection matrix. The transformed objectives may be viewed as a {\em submatrix volume maximization} problem, which involves choosing the optimal $r$-row selection of $\bC\bTi_r$ with the largest possible determinant. Finding this optimum is an NP-hard, intractable combinatorial search over all possible $r$-row submatrices of $\bC\bTi_r$. However, it can be optimized greedily and efficiently via one-time matrix QR factorization requiring $\mathcal{O}(pr^2)\mbox{ and } O(qr^2)$ operations, as described next.

\subsection{ QR pivoting algorithm}
The QR factorization with column pivoting is a greedy submatrix volume optimization scheme that we will use to construct $\bC$ and $\bB$, given $\bTi_r$ and $\bT_r$. The pivoted QR factors any input matrix $\mathbf{V}\in\reals^{r\times p}$ into a unitary matrix $\bQ$, and upper-triangular matrix $\bR$, and column permutation matrix $\mathbf{P}$ so that the permuted matrix $\mathbf{V P}$ is better conditioned than $\mathbf{V}$
\begin{equation}
\label{eqn:qr}
\mathbf{VP} = \bQ\bR.
\end{equation}
\rev{This factorization provides a numerically stable way to compute the determinant as the product of diagonal entries in $\bR$}{However, we seek a well-conditioned row permutation of $\bC\bTi_r$. Consider the input $\mathbf{V}=(\bC\bTi_r)^*$ to the QR factorization, and the leading $r\times r$ square submatrices of the permuted input on both sides of~\eqref{eqn:qr},$\mathbf{\hat{V}}_P$ and $\mathbf{T}$ 
\begin{equation}
\begin{bmatrix}
\mathbf{\hat{V}}_P\mid & \ast~
\end{bmatrix} = [\mathbf{Q}] 
\begin{bmatrix}
\mathbf{T}\mid & \ast~
\end{bmatrix}.
\end{equation}
}
Each iteration of pivoting works by applying orthogonal projections to successive columns of $\mathbf{V}$ to introduce subdiagonal zeros in $\bR$. For our purposes, $\mathbf{P}$ plays the crucial role: at each step $\mathbf{P}$ stores the column ``pivot'' index of the column selected at each iteration to guarantee the following {\em diagonally dominant} structure in $\bR$
	\begin{equation}
	|R_{ii}|^2 \ge \sum_{j=i}^k |R_{jk}|^2; \quad 1\le i \le k \le p.
	\end{equation}
\rev{}{
Observe that the quantity of interest, the determinant of the row-selected submatrix $\mathbf{\hat{V}}_P$ corresponding to the subset selection of measurements, now satisfies
\begin{equation}
|\mathbf{\hat{V}}_P| = |\bQ||\mathbf{T}|   = \prod_{i=1}^r |T_{ii}|,
\end{equation}
since $\bQ$ is unitary and $\mathbf{T}$ is upper-triangular.
}
Because the determinant is the product of these diagonal entries, it can be seen that diagonal dominance guaranteed by the pivoting implicitly optimizes the desired submatrix determinant. 
\rev{To see this, consider the following pivoted QR factorization
\begin{equation}\label{eqn:qr}
\bTi_r^*\mathbf{P} = \bQ\bR.
\end{equation}
Here, column pivoting of the transpose (row selection of direct modes) corresponds to point sensor selection from the ambient high-dimensional space.}{} Thus $\SubC$ is constructed from the first $r$ columns of $\mathbf{P}$ transposed
\begin{equation}\label{eqn:qrC}
\SubC \triangleq (\mathbf{P}_{.,j})^T, \mbox{ where } j:1\rightarrow r.
\end{equation}
\rev{
\begin{algorithm}[t]
	\caption{\small Businger-Golub QR column pivoting of $\bU\in\reals^{r\times n}$. }\label{alg:qrpivot}
	\begin{algorithmic}[1]		\small
	 	\renewcommand{\algorithmicrequire}{\textbf{Input:}}
	 	\renewcommand{\algorithmicensure}{\textbf{Output:}}
	 	\renewcommand{\algorithmiccomment}[1]{$\triangleright$ \textit{#1}}
	 	\REQUIRE $\bU,~p$
	 	\ENSURE  $\Ind$
		\STATE $\Ind \gets [~~]$
		\FOR {$k=1,\dots,p$}
		\STATE $\ind_k = \argmax_{j\notin \Ind} \|\mathbf{v}_j\|_2 $ \hfill \COMMENT{$\mathbf{v}_j$ is the $j$th column of $\bU$}
		\STATE Compute Householder $\tilde{\bQ}$ such that $\tilde{\bQ} \cdot \begin{bmatrix} v_{kk}, \dots, v_{rk} \end{bmatrix}^T = \begin{bmatrix}
		\bR_{kk}, 0 , \dots, 0 \end{bmatrix}^T$  
		\STATE $ \bU \gets \mbox{diag}(\mathbf{I}_{k-1},\tilde{\bQ}) \cdot \bU$ 
			\hfill \COMMENT{remove from all columns the orthogonal projection onto $\mathbf{v}_{\ind_k}$}
		\STATE $\Ind \gets [\Ind,~\ind_k]$
		\ENDFOR
		\RETURN $\Ind$
	\end{algorithmic}
\end{algorithm}
}{}
Actuator selection proceeds similarly to construct a submatrix of $r$ columns of $\bB^*\bT_r$ with maximal determinant, using one additional QR factorization
\rev{}{
\begin{equation}\label{eqn:qr}
(\bT_r^*\bB)\tilde{\mathbf{P}} = \tilde\bQ\tilde\bR.
\end{equation}
}%
The solution $\SubB$ is precisely the leading $r$ columns of $\tilde{\mathbf{P}}$, $\SubB \triangleq \tilde{\mathbf{P}}_{.,j}$, \rev{}{ and we denote by
\begin{equation}\label{eqn:qrB}
\hat\bC = \SubC\bC, \quad \hat\bB = \bB\SubB
\end{equation}
the new measurement and actuation operators obtained in this manner.}

\rev{The QR pivoting procedure is summarized in algorithm~\ref{alg:qrpivot}. Here, we use classical Businger-Golub pivoting~\cite{Businger1965nm}, which applies Householder projections of the input matrix to introduce subdiagonal zeros in $\bR$. In practice, more robust and efficient rank-revealing routines are widely implemented in scientific computing libraries. The operation can be accelerated further using randomization to select pivots~\cite{Drmac2016siam,martinsson2017householder,duersch2017randomized}. 
\noindent\textbf{Historical significance.} It is well-known that the pivoted QR factorization of an input matrix, $\bTi_r^*\mathbf{P}$, is numerically well conditioned for Gaussian elimination and solving least-squares problems~\cite{Businger1965nm}. 
The idea of using the pivoting matrix $\bC$ to select interpolation points for POD modes is first proposed in~\cite{Drmac2016siam}, and its use for sensor placement is detailed in~\cite{manohar2017data}. Due to its favorable numerical properties, this permutation matrix is used to approximate optimal interpolation points for polynomial interpolation~\cite{Sommariva2009qr}.

The performance of QR pivot sensors can be analyzed via their ability to estimate the full state from partial measurements.
Consider the case of sensing in isolation, without actuation. 
The full state estimate is defined by a projection onto direct modes that depends on the chosen sensors $\bC$. The resulting partial observations can be expressed in terms of the balanced coordinates $\ba_r$ as}{ 
The QR pivoting routine is a standard tool in scientific computing for matrix decomposition and linear least-squares problems. We use a block accelerated implementation of classical Businger-Golub pivoting~\cite{Businger1965nm} in MATLAB. Recently QR pivoting was used for interpolating nonlinear terms in EIMs~\cite{Drmac2016siam}, which would otherwise require the evaluation of high-dimensional inner products. In this setting, the interpolation point selection operator is analogous to our selection operator $\SubC$ used with pointwise measurements ($\bC=\mathbb{I}$). The algorithm can be analyzed in terms of the error between the full state and the interpolant approximation at QR pivot interpolation points. The interpolation points can now be written  
}
\begin{equation}
\by = \hat\bC\bx \approx \hat\bC\bTi_r\ba_r,
\end{equation}
where $\bTi_r$ are the POD modes of the reduced model, and $\ba_r$ are the modal coefficients. Recovering the state using the interpolant in the POD basis is accomplished with standard least-squares approximation
\rev{The best estimate of the balanced state from observations, in the least squares sense, is given by
\begin{equation}
\hat\ba_r = (\bC\bTi_r)^{-1}\by,
\end{equation}
and the high-dimensional state estimate can be computed as}{}
\begin{equation}\label{eqn:xproj}
\hat{\bx} = \bTi_r (\hat\bC\bTi_r)^{-1}\by = \bTi_r(\hat\bC\bTi_r)^{-1}\hat\bC\bx.
\end{equation}
This can be expressed as a projection $\bP_C \triangleq \bTi_r(\hat\bC\bTi_r)^{-1}\hat\bC$ of the true state $\bx$ into the observable subspace.
\rev{
\begin{equation}\label{eqn:proj}
\bP_C \triangleq \bTi_r(\hat\bC\bTi_r)^{-1}\hat\bC.
\end{equation}%
}{}
As we shall see, the upper bound on the approximation error 
\begin{equation}\label{eqn:xprojerr}
\|\bx - \bTi_r(\hat\bC\bTi_r)^{-1}\hat\bC\bx\|_2
\end{equation}
is given by $\|(\hat\bC\bTi_r)^{-1}\|_2 = 1/|T_{rr}|$. 
\rev{ Controlling the growth of the latter is closely related to maximizing the determinant, In fact, the two objectives are known as {\em A-optimal} and {\em D-optimal} criteria in experiment design.}{The connection between the latter and maximizing the submatrix determinant can be made explicit in terms of the Hankel singular values of $G$.}

\section{Analysis}
\rev{
We now present lower and upper bounds for the approximation error~\eqref{eqn:xprojerr} given our choices of $\hat\bC$~\eqref{eqn:qrC} and $\bB$~\eqref{eqn:qrB}. 
These bounds rely on the fact that in the case of full observation, $\bC=\bI$, the projection~\eqref{eqn:proj} is simply the approximation in balanced coordinates
\begin{align}
\bx_\star \triangleq \bTi_r\bT_r^*\bx,
\end{align}
resulting in near-optimal bounds on the approximation error. 
\begin{lemma}[Antoulas~\cite{Antoulas2005book}]
The $\ell_2$ error between the full state and its balanced approximation is bounded by}{The best approximation to the state in the span of the direct modes is given by $\bx_\star \triangleq \bTi_r\bT_r^*\bx$ in the ideal measurement scenario $\by = \bx$, i.e. $\hat\bC =\mathbb{I}$. Here the approximation is bounded by the well-known balanced truncation error}
\begin{equation}\label{eqn:bal_bound}
 \|\bx-\bx_\star\|_2 \le 2(\sigma_{r+1}+\dots+\sigma_n),
\end{equation}
where $\sigma_k$ are the Hankel singular values, the diagonal entries of the balanced gramian $\boldsymbol{\Sigma}$~\eqref{Eq:BalancedgramiansHSV}.
\rev{An equivalent result holds for the adjoint system, for which the adjoint state projected onto balanced coordinates is given by $\bz_\star \triangleq \bT_r\bTi_r^*\bz$, and
\begin{equation}\label{eqn:bal_bound2}
 \|\bz-\bz_\star\|_2 \le 2(\sigma_{r+1}+\dots+\sigma_n).
\end{equation}}{
The analysis of empirical QR interpolation in the balanced modes begins with an established result for measurements selected using QR, which states that $\|(\hat\bC\bTi_r)^{-1}\|_2$ at most grows as $\sqrt{p}\mathcal{O}(2^r)$.}
%
\begin{lemma}[Drmac \& Gugercin~\cite{Drmac2016siam}]
The spectral norm of $(\mathbb{S}\mathbf{U})^{-1}$ where $\mathbb{S}$ is computed from the QR factorization~\eqref{eqn:qrC} of the full-rank matrix $\mathbf{U}\in\reals^{p\times r}$ is bounded above
\begin{equation}\label{eqn:qrboundC}
\|(\mathbb{S}\mathbf{U})^{-1}\|_2 \le \frac{\sqrt{p-r+1}}{\sigma_{\min}(\mathbf{U})}\frac{\sqrt{4^r+6r-1}}{3}.
\end{equation}
\end{lemma}%
\rev{}{We generalize this result to the setting of arbitrary linear measurements and actuation, by analyzing the residual between the state and its interpolation in balanced coordinates. Note that the residual between the state and its {\em projection} into balanced modes $\mathbf{v} = \bx - \bx_\star$ satisfies
\begin{equation*}
\bP_C\mathbf{v} = \bP_C\bx - \bTi_r(\hat\bC\bTi_r)^{-1}\hat\bC\bTi_r\bT_r^*\bx_\star = \bP_C\bx - \bx_\star.
\end{equation*}
The interpolation error from QR pivot selection satisfies
\begin{align*}
\|\bx-\bP_C\bx\|_2 &= \|(\mathbf{v}+\bx_\star) - (\bP_C\mathbf{v}+\bx_\star)\|_2 = \|(\bI-\bP_C)\mathbf{v}\|_2 \\
&\le \|\bP_C\|_2\|\bx-\bx_\star\|_2 \\
&\le \|\bTi_r\|_2 \|(\hat\bC\bTi_r)^{-1}\|_2\|\bC\|_2 \|\bx-\bx_\star\|_2. 
\end{align*}%
Substituting~\eqref{eqn:bal_bound},\eqref{eqn:qrboundC} above yields the following result.}
\rev{\begin{theorem}
For any $r$-truncated row permutation matrix $\bC$, the projection error~\eqref{eqn:xprojerr} satisfies the following upper bound
\begin{equation}
\|\bx-\bP_C\bx\|_2 \le \|\bTi_r\|_2 \|(\bC\bTi_r)^{-1}\|_2 \|\bx-\bx_\star\|_2.
\end{equation}
\end{theorem}
\begin{IEEEproof}
Define the residual between full state and its projection into balanced coordinates as $\mathbf{v} = \bx - \bx_\star$. Then
\begin{align*}
\bP_C\mathbf{v} &= \bP_C\bx - \bP_C\bx_\star \\
&= \bP_C\bx - \bTi_r(\bC\bTi_r)^{-1}\bC\bTi_r\bT_r^*\bx_\star \\
&= \bP_C\bx - \bx_\star,
\end{align*}
where we use the fact that the orthogonal projection of $\bx_\star$ onto $\mathcal{R}(\bTi_r)$ is $\bx_\star$ again. Then
\begin{align*}
\|\bx-\bP_C\bx\|_2 &= \|(\mathbf{v}+\bx_\star) - (\bP_C\mathbf{v}+\bx_\star)\|_2 \\
&= \|(\bI-\bP_C)\mathbf{v}\|_2 \\
&\le \|\bP_C\|_2\|\bx-\bx_\star\|_2 \\
&\le \|\bTi_r\|_2 \|(\bC\bTi_r)^{-1}\|_2\|\bC\|_2 \|\bx-\bx_\star\|_2 \\
\end{align*}
\end{IEEEproof}
This logic closely follows that of Chaturantabut and Sorensen~\cite{Chaturantabut2010siamjsc}, in which $\bTi_r$ are restricted to be orthogonal.  This property is not true in general for balanced modes, so $\|\bTi_r\|_2$ is not necessarily equal to $1$. }{}

\begin{theorem} 
The approximation error from interpolating QR-selected observations~\eqref{eqn:qrC} in balanced truncated modes is controlled by the discarded Hankel singular values and the norms of the given measurements and direct modes
\begin{align}
\|\bx-\bP_C\bx\|_2 \le \frac{\|\bC\|_2\|\bTi_r\|_2}{\sigma_{\min}(\bC\bTi_r)}\sqrt{p}\mathcal{O}(2^r) \sum_{i=r+1}^n \sigma_i.  
\end{align}
\end{theorem}
\rev{
\begin{IEEEproof}
Combining Lemmas 3 and 4 with Theorem 1 yields
\begin{align*}
\hspace{-.15in}\|\bx-\bP_C\bx\|_2 &\le \|\bTi_r\|_2 \|(\bC\bTi_r)^{-1}\|_2 \left(2\sum_{i=r+1}^n \sigma_i\right) \\
&\le \|\bTi_r\|_2 \frac{\sqrt{n-r+1}}{\sigma_{\min}(\bTi_r)}\frac{\sqrt{4^r+6r-1}}{3}\left(2\sum_{i=r+1}^n \sigma_i\right) \\
\end{align*} 
\vspace{-.3in}
\end{IEEEproof}

\subsection{Actuator placement}
Actuator placement is completely analogous to sensor placement, where we construct $\bB$ from QR pivoting of the adjoint modes $\bT_r$ instead of the direct modes. Thus, we seek to place actuators along the most controllable directions in state space. Importantly, the controllability gramian in the original system is equivalent to the observability gramian of the adjoint system\hspace{-.05in}
\begin{subequations}
\begin{align}
\dot\bz &= -\bA^*\bz - \bC^*\mathbf{\tilde\by} \\
\tilde\bu &= \bB^*\bz.
\end{align}
\end{subequations}
As before, we desire the best actuation matrix $\bB$ in the adjoint system for estimating the adjoint state $\bz$. 
The adjoint system transforms via the same balanced transformation from above
\begin{subequations}
\begin{align}
\dot{\mathbf{b}}_r &= -\bTi_r^*\bA^*\bT_r\mathbf{b}_r  \\
\tilde\bu &= \bB^*\bT_r \mathbf{b}_r.
\end{align}
\end{subequations}
As previously, the adjoint state may be approximated in balanced coordinates using an analogous projection operator
\begin{equation}\label{eqn:xprojB}
\hat\bz = \bT_r(\bB^*\bT_r)^{-1}\bB^*\bz = \bP_B\bz.
\end{equation}
Substituting $\bT_r$ for $\bTi_r$, $\bz$ for $\bx$, and $\bB^*$ for $\bC$ in the proof for Theorem 1, we obtain the following corollary.

\begin{corollary}
For any $r$-truncated row permutation matrix $\bB^*$, the error from the projection~\eqref{eqn:xprojB} satisfies the following upper bound
\begin{equation}
\|\bz-\bP_B\bz\|_2 \le \|\bT_r\|_2 \|(\bB^*\bT_r)^{-1}\|_2 \|\bz-\bz_\star\|_2.
\end{equation}
\end{corollary}}{The term $\|\bC\|_2 = \|\hat\bC\|_2$ results from information loss when $\bC\ne\mathbb{I}$. An analogous result is obtained for actuator selection by considering the dual problem of estimating the adjoint state from actuation matrix $\hat\bB$ - which is now the \emph{measurement} operator of the adjoint system. The resulting projection operator, $\bP_B \triangleq \bT_r(\hat\bB^*\bT_r)^{-1}\hat\bB^*$, now projects on the span of the \emph{adjoint} modes $\bT_r$. Making appropriate substitutions of $\bP_B$ in the above results yields the following. }
\newtheorem{corollary}{Corollary}
\begin{corollary} 
The approximation error from interpolating QR-selected observations~\eqref{eqn:qrB}  of the adjoint state in balanced truncated modes is controlled by the discarded Hankel singular values and the norms of the given actuators and adjoint modes
\begin{equation}
\|\bz-\bP_B\bz\|_2 \le \frac{\|\bT_r\|_2\|\bB\|_2}{\sigma_{\min}(\bT_r^*\bB)}\sqrt{q}\mathcal{O}(2^r) \sum_{i=r+1}^n \sigma_i.
\end{equation}
\end{corollary}%
\rev{\begin{IEEEproof}
This result immediately follows from the previous corollary and upon substituting $\bT_r$ for $\bTi_r$ in Lemma 4.
\end{IEEEproof}

The above results are for the case when the number of sensors or actuators equal the rank of the truncated model ($p=q=r$), but they generalize to oversampling because the same upper bounds still hold for $p>r$ and $q>r$. }{}%
\rev{\subsection{Log determinant objective}}{}
We now relate the approximation error bounds using QR pivot sensors and actuators to the log determinant objectives. 
\begin{theorem}
Given direct modes $\bTi_r$, QR pivot sensors $\hat\bC$ guarantee the following lower bound for the log determinant \rev{}{
\begin{align}
r \log \frac{9\sigma_{\min}^2(\bC\bTi_r)}{(p-r+1)(4^r + 6r-1)} \hspace{-.15em}+ \hspace{-.3em} \sum_{i=1}^r \log \sigma_i \le \log |\hat\bC\hat\bW_c\hat\bC^T|. \nonumber
\end{align}}
\end{theorem}
\begin{IEEEproof}
Noting the relationship between the singular values of a matrix and its QR factorization, we can express $|\hat\bC\bTi_r|$ in terms of the diagonal entries of its $\bR$ factor
\begin{align}
|\hat\bC\bTi_r| &= \prod_{i=1}^r\sigma_i(\hat\bC\bTi_r) =  \prod_{i=1}^r|T_{ii}| \ge  |T_{rr}|^r,
\end{align}
due to nondecreasing $\sigma_i(\hat\bC\bTi_r)$ for increasing $i$. By squaring the inequality and multiplying by $|\bS_r|$ we obtain
\begin{align}
T_{rr}^{2r}\cdot |\bS_r| &\le |\hat\bC\bTi_r|^2 \cdot |\bS_r|  | \hat\bC\bTi_r\bS_r\bTi_r^*\hat\bC^T| = |\hat\bC\hat\bW_c\hat\bC^T|, \nonumber 
\end{align}
where taking logarithms yields
\begin{align}
 r \log T_{rr}^2 + \sum_{i=1}^r \log \sigma_i \le \log |\hat\bC\hat\bW_c\hat\bC^T|. \nonumber
\end{align}
Because $\|(\hat\bC\bTi_r)^{-1}\|_2 = 1/|T_{rr}|$, the upper bound~\eqref{eqn:qrboundC} in Lemma 2 is the inverse lower bound for $|T_{rr}|$, which can now be substituted above to obtain the final result.
\end{IEEEproof}
An analogous lower bound can be obtained for the objective using QR pivot actuators by appropriately substituting $\hat\bB,\tilde\bR$ and adjoint modes $\bT_r$ in the above proof. 
\begin{corollary}
Given adjoint modes $\bT_r$, $\hat\bB$ satisfies the following lower bound for the log determinant \rev{}{
\begin{align}
r \log \frac{9\sigma_{\min}^2(\bT_r^*\bB)}{(q-r+1)(4^r + 6r-1)} \hspace{-.15em}+ \hspace{-.3em} \sum_{i=1}^r \log \sigma_i \le \log |\hat\bB^T\hat\bW_o\hat\bB|. \nonumber
\end{align} }
\end{corollary}

\section{Results}

We evaluate the selection algorithm in two settings. The first compares QR pivot selections with all possible sensor subset selections in a random state-space model of tractable size. Next we consider an application to closed-loop flow control using LQG control to stabilize unstable Ginzburg-Landau dynamics. The LQG controller with full actuation and sensing is also tractable, and we approximate the $H_2$ optimal placements computed using gradient descent~\cite{Chen:2011} with our QR scheme. 

\subsection{Discrete random state space}

\rev{For our first example, we investigate sparse sensor and actuator placement on random discrete-time state-space systems, generated using the Matlab command \verb|drss|.  
First, we empirically compare the results of QR sensor placement against a brute-force search across all possible sensor placements on a system with $n=25$ states, $p=7$ point sensors, and full-state actuation $\bB = \mathbb{I}$.}{Our first example investigates sensor and actuator selection for random state-space systems with randomized $\bA,\bB,\bC$. First, we compare the results of QR sensor placement against a brute-force search across all possible sensor selections for a system with $n=25$ states and $r=7$ randomized measurements.}
The log determinant objective~\eqref{eqn:logdet_full} is evaluated for all possible choices of 7 sensors, since the system is small enough to explicitly compute the full gramian for all ${n \choose r} = 480,700$ choices of $\hat\bC$. 
These results are binned in Fig.~\ref{fig:brute_force}, and compared with the value resulting from our method (red line). The input to the QR scheme, the balancing modes, are computed only once from the full system. The sensors resulting from our method are observed to be near optimal for the log determinant, exceeding $99.99\%$ of all others\rev{}{, and also good substitutes for $H_2$ optimal sensors.
On average, our method surpasses $99.8\%$ of possible outcomes with a standard deviation of $0.85\%$, over a randomly generated ensemble of 500 model realizations.} 
Therefore, QR sensors are closer to optimal than the analysis suggests.

\begin{figure}[t!]
\centering
\begin{overpic}[width=.48\textwidth]{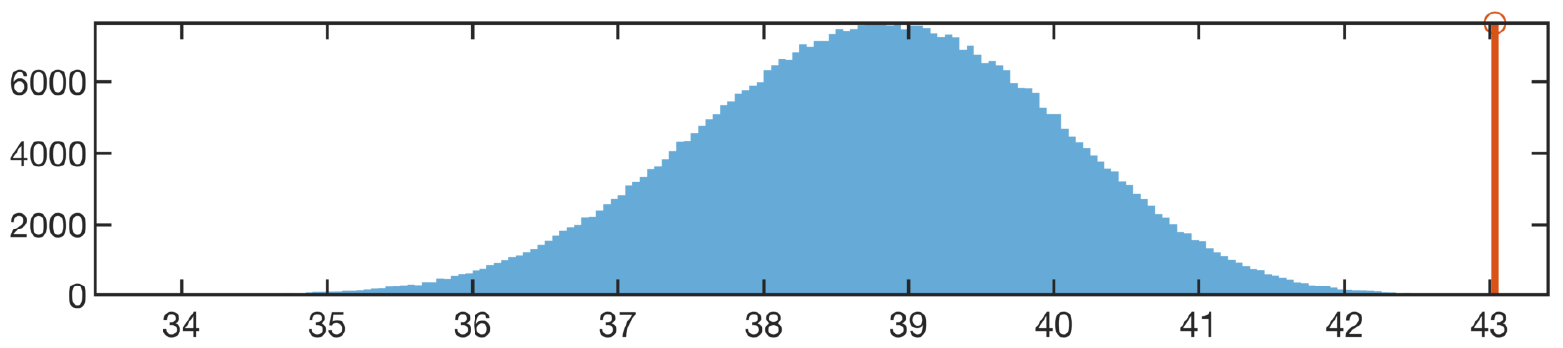}
\put(10,15){$\log | \hat\bC\bW_c\hat\bC^*|$}
\end{overpic}
\begin{overpic}[width=.48\textwidth]{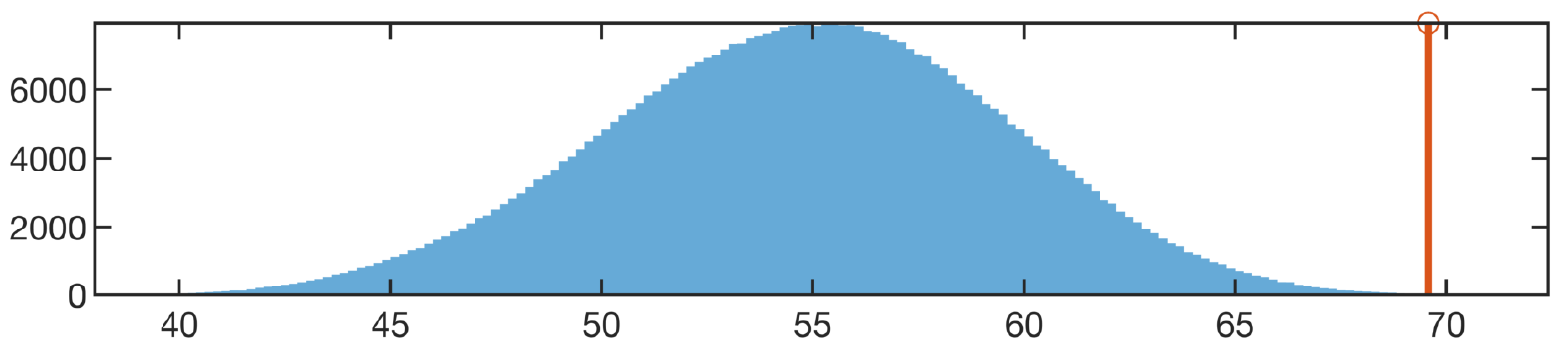}
\put(10,15){$\tr\hat\bC\bW_c\hat\bC^*$}
\end{overpic}
\caption{\color{black} QR pivot sensors (red) greedily maximize the log determinant objective and $H_2$ norms (trace) over all possible selections of 7 sensors out of 25 (blue).} \label{fig:brute_force}
\vspace{-.25in}
\end{figure}

We now investigate performance on a larger random state-space model with $n=100$ states, and likewise initialize the model with \rev{full}{randomized} actuation and sensing such that $p=q=100$.\rev{ with $\bB=\bC=\bI$. }{}
Figure~\ref{fig:rand_plmt} shows the log determinant objective that is being optimized for various sensor and actuator configurations. 
The log determinant of the gramian volume is plotted for the truncated model with QR-optimized sensor and actuator configurations (red circles) and with random configurations (blue violin plots). 
\rev{
In panels (a) and (b), the number of balanced modes is fixed at $r=5$, and the numbers of sensors and actuators are varied, using an extension of QR pivoting for the case when $p>r$ or $q>r$ first described in~\cite{manohar2017data}. }{The truncation level $r$ for the balanced truncation is chosen to match the sensor and actuator budget on the $x$-axis.}
The QR-optimized configurations dramatically outperform random configurations.  \rev{This is less dramatic in the case of optimized partial sensing with full actuation (b), although the QR optimized configurations still outperform the random distribution.  
When the number of sensors and actuators are fixed at $p=q=10$ in (c) }{}As more modes are retained, the chosen sensors and actuators better characterize the input--output dynamics, and their performance gap over random placement increases over all random ensembles, giving empirical validation of our approach.

\begin{figure}[!t]
\centering
{\includegraphics[width=.45\textwidth]{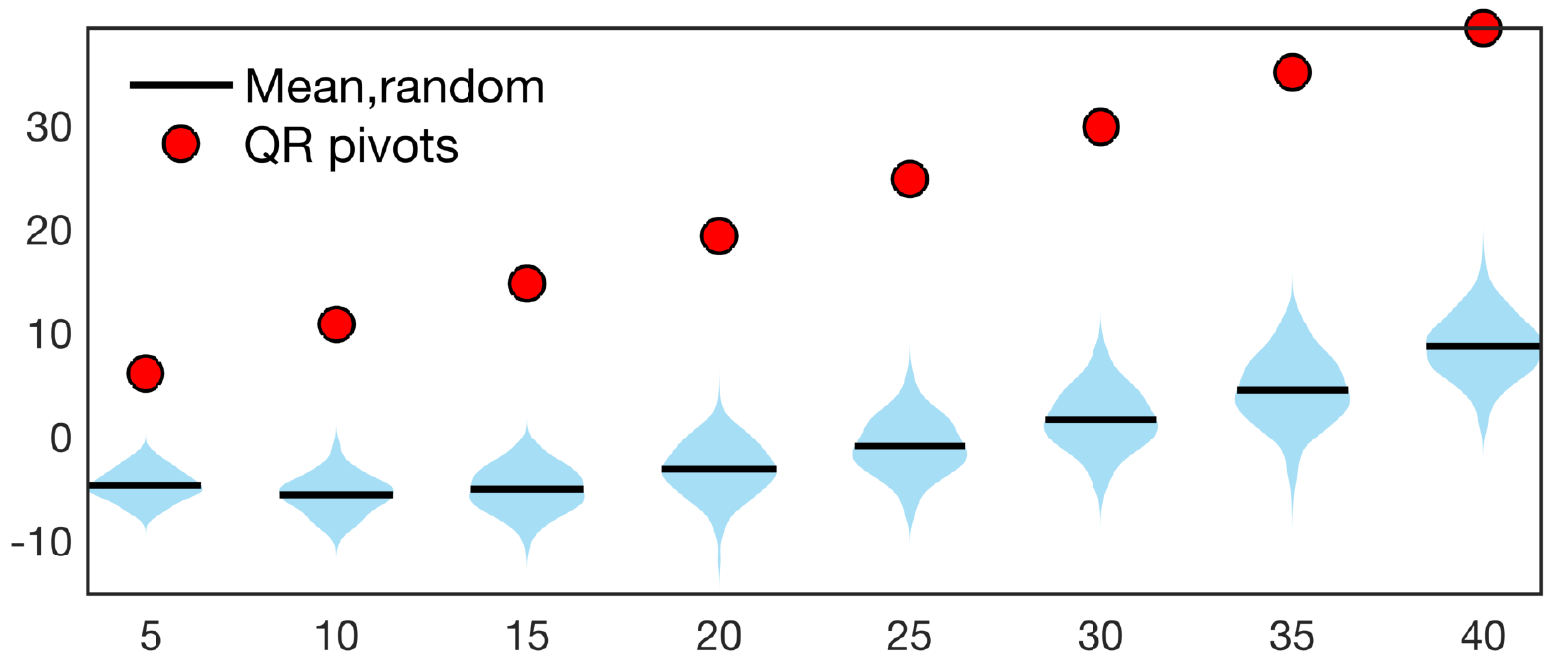}\label{fig_act}}
\caption{Sensor and actuator placement in a random state-space system. The log determinant objective is plotted for  QR-optimized sensor-actuator selections (red) and an ensemble of 200 random sensor-actuator selections (blue violin plots). The truncation level $r$ (also the sensor/actuator budget) varies on the horizontal axis.} \label{fig:rand_plmt}
\end{figure}

Because the system is randomly generated and the dynamics do not evolve according to broad, non-localized features in state-space, many sensors and actuators are required to characterize the system. 
In particular, this is reflected in the slow decay of Hankel singular values. 
By contrast, the next example is generated by a physical fluid flow model, and has coherent structure that allow for a more physical interpretation of sensor and actuator placements with enhanced sparsity.

\subsection{Linearized Ginzburg-Landau with stochastic disturbances}
\rev{As a more sophisticated example, we consider the nonlinear Ginzburg-Landau equation, which models velocity perturbations in a given flow configuration. 
In the case of small perturbations, the flow is well-described by linearized equations
\begin{subequations}\label{eqn:lgle}
\begin{align}
\dot\bx  &= \bA\bx + \bB\bu + \mathbf{D}^{1/2}\mathbf{d} \\
\by &= \bC\bx + \mathbf{N}^{1/2}\mathbf{n},
\end{align}
\end{subequations}
where $\mathbf{d}$ is a stochastic process disturbance and $\mathbf{n}$ is noise present everywhere in the domain. The $H_2$ optimal sensor and actuator placement for this system has been determined by Chen and Rowley~\cite{Chen:2011} using gradient minimization, and hence provides a benchmark comparison for our greedy placement. 

The observability and controllability gramians are not defined for this system since it is unstable. It can, however, be stabilized using linear quadratic Gaussian (LQG) control, which seeks to stabilize the system with minimal input using the cost function
\begin{equation}
\mathbf{J} = \bx^*\hat\bQ\bx + \bu^*\hat\bR\bu,
\end{equation}
where $\hat\bQ$ and $\hat\bR$ refer to user-specified state and input weighting matrices, not to be confused with the matrices in the QR factorization.
The resulting system is stable, therefore we can compute the balancing transformations using Matlab's \verb|balreal| command. 
Specifically, the new system matrices $\bA_K,\bB_K,\bC_K$  are computed using Matlab's \verb|lqg| routine, and its balanced modes are computed using \verb|balreal|. 
Explicit expressions for the controller system matrices are given in~\cite{Chen:2011}. 
In their work, sensing and actuation are defined on smooth spatial Gaussian kernels centered at the sensor/actuator location, and the small-width limit of the kernels corresponds to point sensing and actuation. }{}
%
{\color{black} We consider the closed-loop linearized Ginzburg-Landau model evolving velocity perturbations in a flow, given a controller with full actuation and sensing, which is often not feasible in practice. The equations modeling the plant dynamics are unstable because the system matrix has eigenvalues in the right half plane. The dynamical system matrix $\bA$ is formed from Hermite pseudospectral discretization of the linearized Ginzburg-Landau operator
\begin{equation}
A \triangleq -\nu\frac{\partial}{\partial \xi} + \mu(\xi) + \beta \frac{\partial}{\partial \xi^2}.
\end{equation}
The spatial grid $\boldsymbol{\xi}\in\reals^n$ is discretized at the $n=100$ roots of Hermite polynomials, and $\nu,\beta,\mbox{ and } \mu(\xi)$ are advection, diffusion and wave amplification parameters. Each $i$th sensor $\xi_s$ (row of $\bC_2$) and actuator at $\xi_a$ (column of $\bB_2$) are weighted by Gaussian kernels and the trapezoidal integration weights $\mathbf{M}$ 
\begin{equation}
\bC_{2_i} \triangleq \left[e^{-\frac{(\boldsymbol{\xi} -\xi_s)^2}{\sqrt{2}\sigma}}\right]^T\mathbf{M}, \quad \bB_{2_i}\triangleq e^{-\frac{(\boldsymbol{\xi} -\xi_a)^2}{\sqrt{2}\sigma}}.
\end{equation}
The linear quadratic Gaussian (LQG) controller stabilizes the dynamics by minimizing the $H_2$ optimal cost function $J(\bx,\bu) = \bx^T\hat\bQ\bx + \bu^T\hat\bR\bu$, where $\hat\bQ$ and $\hat\bR$ are user-specified weight matrices. The output $\bu$ of the LQG controller, given by
\begin{equation}\label{eqn:lqg}
\begin{bmatrix}\dot{\hat{\bx}} \\ \bu\end{bmatrix} =  \begin{bmatrix}\bA-\bB_2\mathbf{F}-\mathbf{L}\bC_2 & \mathbf{L} \\ -\mathbf{F} & 0 \end{bmatrix} \begin{bmatrix}\hat{\bx} \\ \by\end{bmatrix},
\end{equation}
stabilizes the dynamics given white noise stochastic disturbance $\mathbf{d}$ and noise $\mathbf{n}$ at all sensors and actuators, with covariances $\mathbf{V} = 4\cdot 10^{-8}\mathbb{I}$ and $\mathbf{W} = \mathbb{I}$.
Since every state is observed and actuated, each spatial gridpoint corresponds to one $\xi_a$ and $\xi_s$. The idea is to preserve as much of this ``ideal'' controller as possible using a subset of the original sensors and actuators
\begin{subequations}\label{eqn:plant}
\begin{align}
	\dot\bx &= \bA\bx + \bB_2\SubB\bu  + \mathbf{W}^{1/2}\mathbf{d}  \\ 
	\by &= \SubC\bC_2\bx + \mathbf{V}^{1/2} \mathbf{n},
\end{align}
\end{subequations}
which is similar to our original problem formulation. Hence we can perform balanced model reduction on either impulse responses or the controller directly, and then QR pivoting to optimize placements. In this formulation, $\bu$ is the output and $\by$ is the input, which encapsulates the notion of maximizing the gain from feedback to $\bu$ to stabilize the dynamics, which translates to the observabillity of~\eqref{eqn:lqg}. Thus we compute gramians and adjoint, direct modes of the LQG matrices $\bA \triangleq \bA - \bB_2\mathbf{F} - \mathbf{LC}_2, \bB \triangleq \mathbf{L},  \bC \triangleq -\mathbf{F}$.

We compare our approach to established gradient descent techniques for computing the $H_2$ optimal controller and sensor-actuator placements simultaneously. The particular algorithm for comparison is the optimal placement for this model determined using the gradient descent scheme of Chen and Rowley~\cite{Chen:2011}.%
}

\begin{figure}[t!]
\centering
\begin{overpic}[width=.45\textwidth]{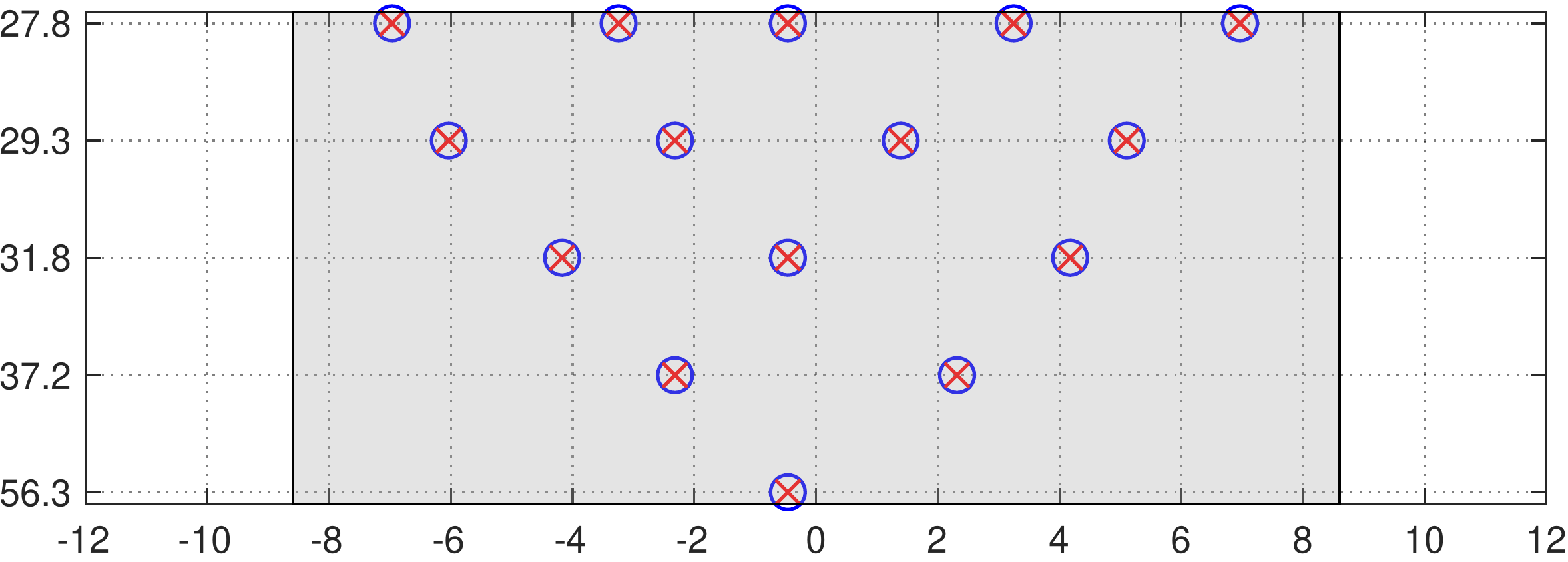}\put(6,30){(a)QR}\end{overpic}
\label{fig_first_case}\\\vspace{-1mm}
\begin{overpic}[width=.45\textwidth]{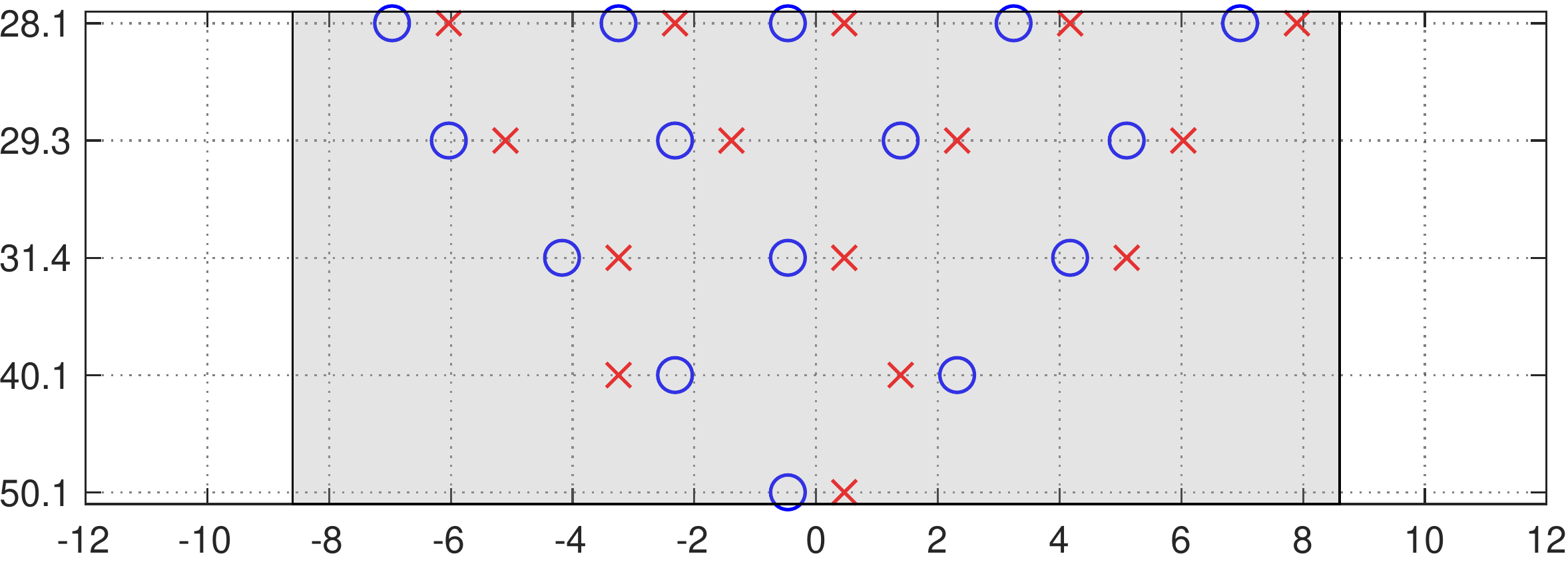}\put(6,30){(b)QR}\end{overpic}
\label{fig_second_case}\\\vspace{-1mm}
\begin{overpic}[width=.45\textwidth]{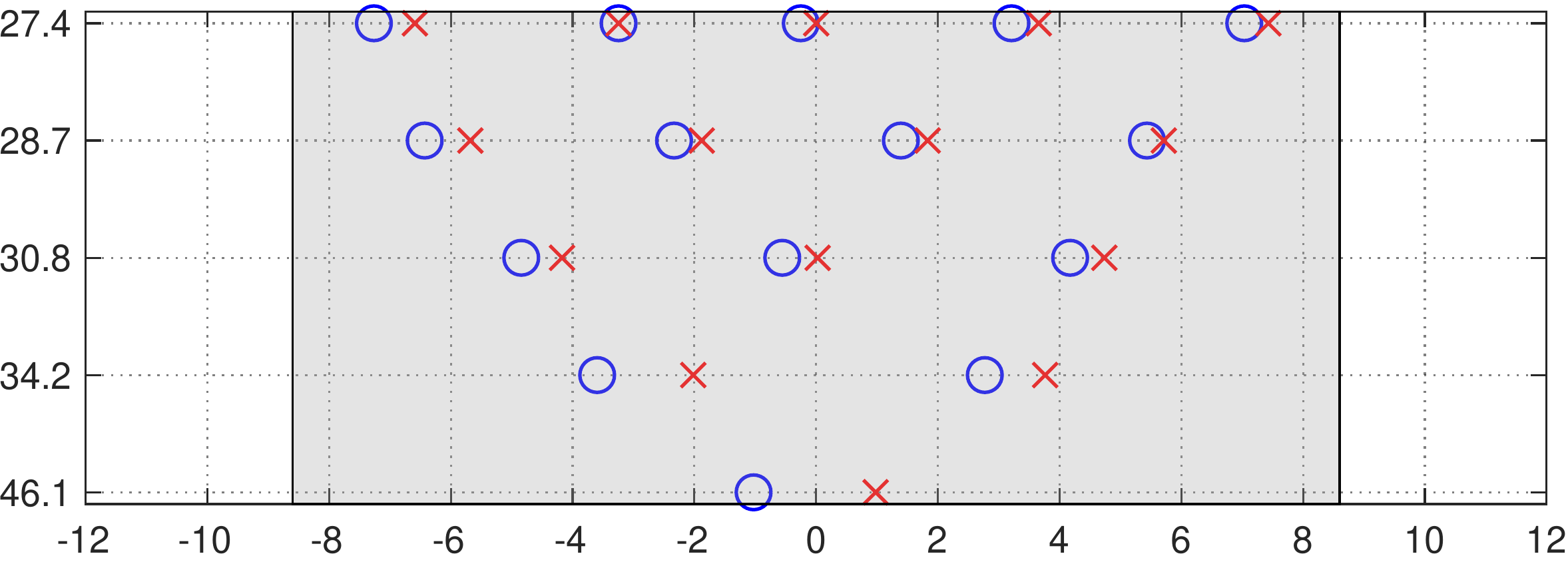}\put(7,30){(c)$H_2$ opt}\end{overpic}
\label{fig_third_case}
\caption{Sensor ({\color{red}$\times$}) and actuator ({\color{blue}$\circ$}) placement for linearized Ginzburg-Landau. 
Each row corresponds to the optimized placement for budgets of 1-5 sensors and actuators. 
Placements based on QR pivoting of balanced truncated modes~(a) closely approximate the $H_2$ norms of the placements determined using gradient descent~(c).  The QR method can be modified to place sensors and actuators to avoid collocation (b). } \label{fig:lgle_sens}
\end{figure}%

Their $H_2$ norm optimization scheme permits placement of sensors and actuators at locations that may not be grid points. The major drawback is that each Newton iteration requires solving $2r$ $n\times n$ Lyapunov equations until convergence\rev{}{, although recent work simplifies this to 2 equations per iteration~\cite{colburn2011gradient}.} Furthermore, the procedure requires an ensemble of random initial conditions to avoid converging to a local minimum. In~\cite{Chen:2011}, the optimal placement is computed using  conjugate gradient optimization for the same spatial discretization $n=100$, which becomes computationally expensive as the grid resolution increases.
\rev{}{%
In this case, gradient descent is more costly than balancing the fully actuated and observed system, which comes at a one-time cost of solving 2 Lyapunov equations for the gramians, and 2 Riccati equations for the LQG gain matrices ($O(n^3)$ each). Therefore, our algorithm is  sensible when the grid discretization is sufficiently fine. Furthermore, our solution is a good starting point for the convergence of the gradient descent scheme, thus eliminating the need for optimization over a large ensemble of randomized starting points. 
QR pivoting runtime scales as $\mathcal{O}(nr^2)$ and the deviation of the resulting placement from the $H_2$ optimum (fig.~\ref{fig:lgle_sens}) decreases with increasing $r$. }
\rev{The structure of the controller is shown in Fig.~\ref{FIG:plant}. Here $\mathbf{j}$ and $\mathbf{w}$ are concatenated vectors containing terms in the cost function, and the concatenated sources of  process and measurement disturbance, as defined below
\begin{equation*}
\mathbf{j} = \begin{bmatrix} \hat\bQ^{1/2}\bx \\ \hat\bR^{1/2}\bu \end{bmatrix}, \mathbf{w} = \begin{bmatrix} \mathbf{d} \\ \mathbf{n} \end{bmatrix}.
\end{equation*}
By minimizing $\mathbf{J}$, the controller minimizes the gain from the disturbances to the weighted cost vectors.
\begin{figure}[t!]
\centering
\begin{tikzpicture}[node distance=2cm]
\tikzstyle{block} = [rectangle, minimum width=2cm, minimum height=1cm,text centered, draw=black]
\tikzstyle{approx} = [rectangle,  minimum width=2cm, minimum height=1cm,text centered, draw=black]
\tikzstyle{arrow} = [->,>=stealth]
\node (plant) [block] {Plant $P$};
\node (klqg) [block, below of=plant] {$K$};
\node [coordinate,left of=plant,yshift=-.15cm] (lbend) {};
\node [coordinate,right of=plant,yshift=-.15cm] (rbend) {};
\node [coordinate,right of=plant,yshift=.15cm] (r) {$w$};
\node [coordinate,left of=plant,yshift=.15cm] (l) {};
\draw [arrow] ([yshift=-.15cm]plant.east) -- node[anchor=north]{$\by$}(rbend) |- (klqg);
\draw [arrow] (klqg) -|  (lbend) -- node[anchor=north]{$\bu$}([yshift=-.15cm]plant.west);
\draw [arrow] ([yshift=.15cm]plant.east) -- node[anchor=south]{$\mathbf{j}$}(r);
\draw [arrow] (l) -- node[anchor=south]{$\mathbf{w}$} ([yshift=.15cm]plant.west);
\end{tikzpicture}
\caption{LQG controller for linearized Ginzburg-Landau\label{FIG:plant}}
\end{figure}}{}%
\begin{figure}[!t]
\centering
\hspace{-.65in}
{\begin{overpic}[width=.35\textwidth]{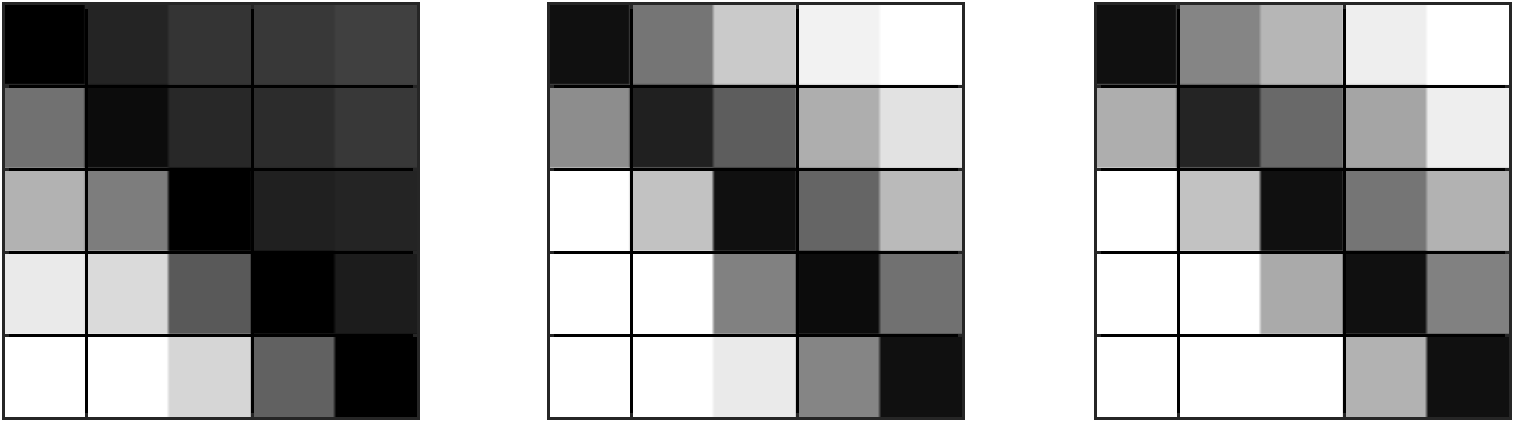}
\put(7,29) {\small$\omega = 10^{-1}$}
\put(45,29) {\small$\omega = 10^{1}$}
\put(80,29) {\small$\omega = 10^{3}$}
\put(-15,29) {(a)$H_2$opt}
\label{fig:lqg_gain_opt}\end{overpic}}
\hspace{-.45in}
\llap{\shortstack{%
        \includegraphics[scale=.325]{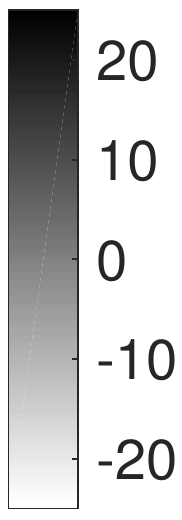}\\
        \rule{0ex}{0in}%
      }
  \rule{-.75in}{0ex}}\\
  \hspace{-.65in}
{\begin{overpic}[width=.35\textwidth]{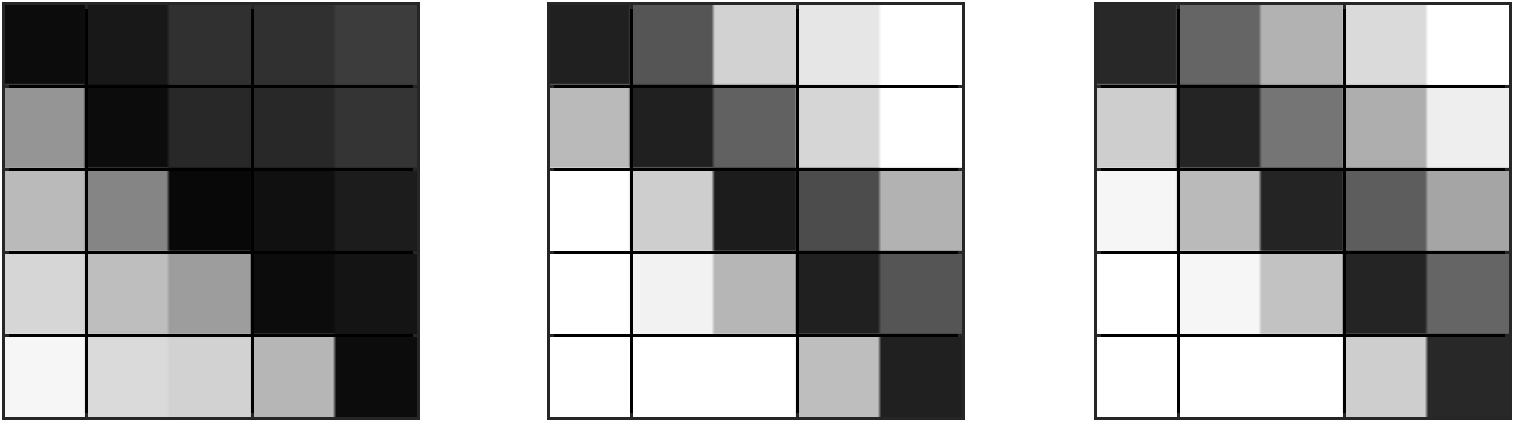}
\end{overpic}\put(-210,29) {(b) QR}
\label{fig:lqg_gain_qr}}
\hspace{-.45in}
\llap{\shortstack{%
        \includegraphics[scale=.325]{lqg_gain_colorbar}\\
        \rule{0ex}{0in}%
      }
  \rule{-.75in}{0ex}}
\caption[LQG gain for five-input, five-output system.]{LQG gain (dB) for a system with 5 sensors and actuators. Each block shows the gain from a signal $\exp(i\omega t)$ in sensor $k$ (column) to actuator $j$ (row), ordered upstream to downstream.} \label{fig:lqg_gain}
\vspace{-.2in}
\end{figure}

Figure~\ref{fig:lgle_sens} plots sensor and actuator configurations from the QR algorithm and $H_2$ gradient optimization, which are compared with the $H_2$ optimal placements in \cite{Chen:2011}.
The resulting placements for the cases $r=1$ to $r=5$ sensors and actuators are plotted vertically, and the horizontal axis is the spatial domain $\xi\in[-12,12]$ with a shaded wave amplification region \rev{}{ in which fluid perturbations are amplified}. For each value of $r$, we apply QR pivoting to the rank $r$ truncated balanced modes. QR pivoting collocates sensors and actuators, indicating that $\bA$ is approximately symmetric and hence the direct and adjoint modes (pictured in Fig.~1) are identical up to a scaling factor. \rev{}{In practice, sensors are often slightly downstream to account for time delays, }so we enforce via the pivoting procedure that sensors are not placed at previously chosen actuators.
The $H_2$ norms of the resulting placement on the $y$-axis indicate that the QR selections closely approximate the optimal placements.  
The $H_2$ optimal placement~\cite[Fig. 4]{Chen:2011} of five sensors and actuators, with $H_2$ norm 27.4, agrees exactly with the $H_2$ optimum and is closely approximated by the QR pivoted placement (27.8). 
\rev{Our results quantitatively agree with the classical $H_2$ optimal placements. }{} 

Figure~\ref{fig:lqg_gain} compares controller performance between QR pivoting and the $H_2$ optimum via the LQG gain of a given signal from each sensor to each actuator. The LQG gains are identical to those produced by the $H_2$ optimal method of Chen and Rowley~\cite[Fig. 5]{Chen:2011}. The diffusive nature of the dynamics favors nearly collocating the sensors and actuators, since the high-frequency oscillations mostly propagate to the nearest actuator. This confirms that our framework is useful for optimizing sensors and actuators. Balanced truncation applied to the closed loop system is critical to achieving this, since the open loop dynamics are unstable and it is shown in~\cite{Chen:2011} that the dominant eigenmodes of the dynamics lead to vastly suboptimal placements.

\section{Discussion and Outlook}

In this work we develop scalable sensor and actuator \rev{placement}{selection} whose runtime scales  linearly with the number of state variables, after a one-time offline computation of the balanced modes. Our approach relies on balanced model reduction~\cite{Moore1981ieeetac,Willcox2002aiaaj,Rowley2005ijbc}, which hierarchically orders modes by their observability and controllability. We extend EIMs to interpolate the low-rank balancing modes of the system and determine maximally observable and controllable locations (sensor \& actuators) in state space. The performance of this algorithm is demonstrated on random state-space systems, and optimal $H_2$ control of the linearized Ginzburg-Landau model.
Our optimized placements vastly exceed the performance of random placements, and closely approximate $H_2$ optimal placements computed by costly gradient minimization schemes, but achieved at a fraction of the runtime. 



Sensors and actuators are critical for feedback control of large  high-dimensional complex systems.
This work advocates sensor and actuator \rev{placement}{selection} using QR pivots of the direct and adjoint modes of a system's balancing transformation. The resulting placement is empirically shown to preserve the dynamics of the full system. The method has deep connections to system observability, controllability, modal sampling methods and classical experimental design criteria. Furthermore, QR pivoting is more computationally efficient than leading greedy and convex optimization methods, and thus critically enlarges the search space of possible \rev{candidate placements}{selections}. This is particularly valuable in spatiotemporal models where high-resolution grids generate a large number of states, and balanced modes and QR method exploit the spatial structures.

This work opens a variety of future directions in pivoting sensor and actuator optimization. Rapid advances in data collection yield extremely large search spaces, for which the computation of balanced modes and QR pivoting may be accelerated using randomized linear algebra. Our method relies on a known model of the dynamics, but it would also be interesting to generalize the method to data-driven system identification models. In addition, point sensors and actuators are simplifications of constrained or nonlinear sensing and actuation that may occur in practice. Nonlinear sensing constraints remain an open challenge.

\bibliographystyle{IEEEtran}
\bibliography{references}
\end{document}